\documentclass[a4paper, 12pt, fleqn]{article}

\usepackage[fleqn]{amsmath}
\usepackage{alltt}
\usepackage{amssymb}

\newtheorem{thm}{Theorem}
\newtheorem{lemma}[thm]{Lemma}

\begin{document}

\title{Characterization of SDP Designs That Yield
Certain Spin Models}

\author{Carl Bracken and Gary McGuire\\
Department of Mathematics\\
NUI Maynooth\\
Co. Kildare\\
Ireland}

\maketitle
\begin{abstract}
\noindent
We characterize the SDP designs that give rise to
four-weight spin models with two values.
We prove that the only such designs are the symplectic
SDP designs.
The proof involves analysis of the
cardinalities of intersections
of four blocks.
\end{abstract}

\newpage

Keywords: symmetric difference property, design, spin model,
symplectic.

\newpage

\section{Introduction}
The concept of a \emph{spin model} was introduced by Jones \cite{J}.
The concept was generalised to \emph{two-weight} spin models
by Kawagoe, Munemasa and Watatani  \cite{KMW}, 
and further generalised to \emph{four-weight} spin models
by Bannai and Bannai  \cite{BB}.

Guo and Huang  \cite{GH} considered certain types
of four-weight spin models, which they called
``four-weight spin models with exactly two values on $W_2$''
(we will explain this further in section 2).
They showed a connection with symmetric designs.
Bannai and Sawano  \cite{BS} showed that the existence
of a four-weight spin model with exactly two values on $W_2$
is equivalent to the existence of a quasi-3 design
with certain properties (see Theorem \ref{bs1} below),
strengthening a result of Guo-Huang \cite{GH}.
They (implicitly) raised the question of whether SDP designs,
which are known to be quasi-3, would satisfy these
properties.  In this paper we will answer this question.
Our main result is as follows.
\bigskip

{\bf Theorem.}  {\sl An SDP design $D$
satisfies all the conditions of Theorem $1$, and thus
corresponds to a four-weight spin model, if and only if
$D$ is equivalent to the symplectic SDP design.}

\bigskip


In section 2 we give some background and definitions.
In section 3 we will prove some preliminary results 
about SDP designs.
In section 4 we prove that the only SDP designs satisfying
condition 2 of Theorem \ref{bs1} are 
(up to design equivalence)
the symplectic designs.
In section 5 we prove that the symplectic SDP designs
satisfy condition 3 of Theorem \ref{bs1},
and then our main result follows.


\section{Background}

We first give the definition of a four-weight spin model.
Matrices will be indexed by the elements of a finite set $X$.
We use $A\circ B$ to denote the Hadamard product of
matrices $A=(A(\alpha, \beta))_{\alpha,\beta\in X}$ and 
$B=(B(\alpha, \beta))_{\alpha,\beta\in X}$,
which is the matrix whose $(\alpha,\beta)$
entry is equal to $A(\alpha, \beta) B(\alpha, \beta)$.
Let $I$ denote the identity matrix, 
let $J$ denote the all-1 matrix,
and let $A^T$ denote the transpose of $A$.

\bigskip

{\bf Definition.}  Let $X$ be a finite set with $n$ elements,
and let $D$ be a real number satisfying $D^2=n$.
We say that $(X, W_1, W_2, W_3, W_4)$ is a
\emph{four-weight spin model} of size $n$ if
each $W_i$ ($1\leq i \leq 4$) is an $n$-by-$n$
matrix with complex entries and the following conditions hold:
\begin{enumerate}
\item $W_1^T \circ W_3 = W_2^T \circ W_4 = J$
\item $W_1 W_3 = W_2 W_4 = nI$
\item \begin{enumerate}
\item $W_1 Y^{41}_{\alpha \beta}=D W_4(\alpha,\beta)\
Y^{41}_{\alpha \beta}$ for all $\alpha, \beta \in X$
\item $W_1^T Y^{14}_{\beta \alpha}=D W_4(\beta,\alpha)\
Y^{14}_{ \beta \alpha}$ for all $\alpha, \beta \in X$
\end{enumerate}
\end{enumerate}
where $Y^{ij}_{\alpha \beta}$ is the $n$-dimensional
column vector whose entry in position $\gamma$ is
$W_i(\alpha,\gamma)W_j(\gamma,\beta)$.

\bigskip

The apparent lack of symmetry in 3(a) and 3(b) is
explained in \cite{BB}, where they show that each of
these equations is equivalent to seven others.

If $W_1=W_2$ and $W_3=W_4$ then the definition above
reduces to the definition of a two-weight spin model
in \cite{KMW}.
If in addition we assume that $W_1$ and $W_3$ are symmetric,
the definition reduces to the definition of a spin model
in \cite{J}.

The papers \cite{GH} and \cite{BS} consider the case
that $W_2$ has only two distinct entries, with each entry
appearing the same number of times in each row and column.
They refer to this case as a \emph{four-weight spin model
with exactly two values on $W_2$}.
This is the case under consideration in this paper.

\bigskip

A block design with parameters 2-$(v,k,\lambda)$ is
often called a symmetric design if $b=v$.
Following \cite{CvL} we shall use the term
``square'' design for a symmetric design.
We refer the reader to \cite{CvL} for the basic properties of
block designs.
A square design is said to be {\it quasi-$3$} {\it for points}
if the number of blocks incident with three distinct
points takes only two values.
Such designs seem to have been first 
considered in Cameron \cite{C}; see also \cite{BM} for a survey
of quasi-3 designs.   We shall say that a 
square design is {\it quasi-$3$ for blocks} if the number of
points in the intersection of any three distinct blocks takes only
two values. 
A design is quasi-3 for blocks if and only if 
the dual design is quasi-3 for points.

We index the rows of an incidence matrix of a design
by the blocks, and the columns by the points.
When we speak of the \emph{sum} of a number of blocks, we mean the
sum modulo 2 of the rows of the incidence matrix
corresponding to those blocks.

An SDP (symmetric difference property)
design is a square $(v,k,\lambda)$ design with
the property that the symmetric difference 
(or sum) of any
three blocks is either a block or a block complement.
It follows immediately from the definition that SDP
designs are quasi-3 for blocks.
For, the identity
$$|B+B'+B''|=|B|+|B'|+|B''|-2|B\cap B'|-2|B\cap B''|
-2|B'\cap B''|$$
$$\qquad \qquad\qquad  +4|B\cap B' \cap B''| \qquad \qquad (\clubsuit)$$
shows that there are only two possibilities for
$|B\cap B' \cap B''|$.
It is shown in \cite{BM} that all SDP designs are also quasi-3
for points.

We  now state the result of \cite{BS} on the case
that $W_2$ has only two distinct entries.


\begin{thm} \label{bs1} $\cite{BS}$
Let $W_2=\alpha A +\beta (J-A)$, where $\alpha, \beta$
are distinct nonzero complex numbers, and $A$ is a
$(0,1)$-matrix with the property that each row and 
column has exactly $k$ ones, where $2\leq k \leq n-2$.
Let $X$ be a finite set with $n$ elements,
and let $D$ be a real number satisfying $D^2=n$.
Then $W_2$ defines a four-weight spin model if and only
if $A$ is the incidence matrix of a square $(n,k,\lambda)$
design $D(X,{\cal B})$ which satisfies the following
three properties:
\begin{enumerate}
\item $D(X,{\cal B})$ is quasi-3 for blocks with triple
intersection sizes
$$x=\frac{k\lambda + \lambda - (k-\lambda)\sqrt{k-\lambda}}{n},\ \
y=\frac{k\lambda + \lambda + (k-\lambda)\sqrt{k-\lambda}}{n}.$$
\item For any set ${\cal S}\subseteq {\cal B}$ 
of four blocks, an even number of the
four $3$-subsets of ${\cal S}$ have triple intersection size $x$.
\item There exists a 1-1 correspondence 
$\phi : X \longrightarrow {\cal B}$
with the property that for any three points
$a,b,c \in X$,
the number of blocks containing $\{a,b,c\}$ is
$|\phi(a)\cap \phi(b)\cap \phi(c)|$.
\end{enumerate}
Moreover, if conditions 1, 2 and 3 hold, then
$\alpha, \beta$ and $W_1$ are determined by $D$ and $k$.
In particular, $\alpha=-\beta$ if and only if
$n=4q^2$ where $q$ is an even integer.
\end{thm}

Guo and Huang \cite{GH} point out that the $(16,6,2)$
SDP design satisfies conditions 1, 2 and 3 of
Theorem \ref{bs1}, and thus gives an example of a
four-weight spin model
with exactly two values on $W_2$.
Bannai and Sawano \cite{BS} 
showed that the other (non-SDP) $(16,6,2)$ designs
do not satisfy conditions 1, 2 and 3 of Theorem \ref{bs1}.
They state that 
it is known that SDP designs are quasi-3 for blocks
(satisfy condition 1 of Theorem \ref{bs1}),
but as we said above this follows from the definition.
They appear to be wondering
whether all SDP designs satisfy
conditions 1, 2 and 3 of Theorem \ref{bs1}.
We investigate this question in this paper.
We will show that, although the number of 
nonisomorphic SDP designs
grows exponentially with $m$,
there is  one and only one SDP design (up to isomorphism)
satisfying the three conditions.

\section{On SDP Designs}

It was shown by Kantor \cite{K1} that any SDP design must
have parameters
$$(v,k,\lambda)=
(2^{2m},2^{2m-1}-2^{m-1},2^{2m-2}-2^{m-1}).$$ 
There is one particular SDP design of interest to us, 
which is called the \emph{symplectic} SDP design.
It is constructed using the $2^{2m}$ quadratic forms
that polarise to a given nondegenerate symplectic
bilinear form on a $2m$-dimensional vector space over $GF(2)$,
see \cite{CS}, \cite{CvL} or \cite{K1}.
The symplectic design has a 2-transitive
automorphism group.
Kantor \cite{K2} showed that
the number of nonisomorphic SDP designs
grows exponentially with $m$.

We recall that a regular
Hadamard matrix is a Hadamard matrix with
constant rowsums.  Such a matrix of size $4u^2$ gives rise
to a square
2-$(4u^2,2u^2-u,u^2-u)$ design (replacing $-1$ by $0$
and perhaps complementing).
When $u=2^{m-1}$ these parameters are the same as
the SDP parameters.
Taking Kronecker products of the 4-by-4 matrix $J-2I$
with itself results in the symplectic SDP designs
(this description is due to Block \cite{B}).

\begin{thm}\label{sumoffour}
Let $D(X,{\cal B})$  be the $2$-$(4u^2,2u^2-u,u^2-u)$ design induced
by a regular Hadamard matrix of size $4u^2$.
Then $D$ is an SDP design if and only if the
sum of any four blocks is a vector of weight
$0$, $2u^2$, or $4u^2$.
\end{thm}

Proof:
Let $v=4u^2$.
First suppose $D(X,{\cal B})$  is an SDP design.
Let $B_1, B_2, B_3\in {\cal B}$ 
be three distinct blocks of $D$.
Then $B_1+B_2+B_3=B$ or $B+j$, where $j$ denotes
the all-1 vector and $B\in {\cal B}$.
If $B_4$ is a block not equal to any of $B_1, B_2, B_3$, then
$$B_1+B_2+B_3+B_4=B+B_4\ \ \textrm{ or }\ \  B+B_4+j.$$
Since the sum of any two distinct blocks has weight $v/2$,
the weight of $B+B_4$ (and $B+B_4+j$) will be one
of $0$, $v/2$, or $v$ (as $B$ could equal $B_4$).

Conversely, suppose $D(X,{\cal B})$  is not an SDP design.
Then there exist $B_1, B_2, B_3\in {\cal B}$
such that $B_1+B_2+B_3$ is not a block or a
block complement.
Let
$${\cal H} = \{B_2+B_3\} \cup \{B_1+B_i: 2\leq i \leq v\}$$
and let 
$$\overline{{\cal H}}=\{x+j:x\in {\cal H}\}.$$
Then $|{\cal H}|=|\overline{{\cal H}}|=v$,
and ${\cal H}\cap \overline{{\cal H}}=\emptyset$.

For the sake of contradiction,
assume that the sum of any four blocks of $D$ has weight
$0$, $v/2$, or $v$. 
Then the sum of any two distinct
elements of ${\cal H}$ has weight $v/2$,
and the same applies to any two distinct elements of 
$\overline{{\cal H}}$.
Also, if $x\in {\cal H}$ and $y+j\in  \overline{{\cal H}}$,
then $x+y+j$ has weight $v/2$ (unless $x=y$
in which case $x+y+j=j$ has weight $v$).

It follows that any two distinct elements of
${\cal H}\cup \overline{{\cal H}}$ have 
Hamming distance at least $v/2$.
Adding the all-0 vector to these vectors yields
a binary $(v,2v+1,v/2)$ code, which violates
the Plotkin bound (see \cite{MS} chapter 2).
This contradiction completes the proof.

\hfill $\Box$

\bigskip

We shall use Theorem \ref{sumoffour} to calculate
the possible quadruple intersection sizes of blocks
in an SDP design.

\begin{thm}\label{quadint}
Let $D(X,{\cal B})$  be a $2$-$(4u^2,2u^2-u,u^2-u)$ SDP design,
where $u=2^{m-1}$.
Then the cardinality of the intersection of four distinct
blocks takes one of the following seven values:
\begin{enumerate}
\item 0
\item $u^2/2 - u$
\item $u^2/4$
\item $u^2/4-u/4$
\item $u^2/4-u/2$
\item $u^2/4-3u/4$
\item $u^2/4-u$.
\end{enumerate}
\end{thm}

Proof:  
Let $B_1, B_2, B_3, B_4\in {\cal B}$ be four distinct
blocks of $D$.  Let
$$\alpha=|B_1\cap B_2\cap B_3|$$
$$\beta=|B_1\cap B_2\cap B_4|$$
$$\gamma=|B_1\cap B_3\cap B_4|$$
$$\delta=|B_2\cap B_3\cap B_4|,$$
and let
$$q=|B_1\cap B_2\cap B_3\cap B_4|.$$
It follows easily from $(\clubsuit)$
that each of $\alpha,\beta,\gamma,\delta$, is equal to either
$x=u^2/2-u$ or $y=u^2/2-u/2$.
Let $w$ be the weight of the vector
$B_1+B_2+B_3+B_4$.
Then (by the obvious generalisation of $(\clubsuit)$ to
four blocks)
$$w=4(2u^2-u)-12(u^2-u)+4(\alpha+\beta+\gamma+\delta)-8q. \quad (\diamondsuit)$$

By Theorem \ref{sumoffour}, $w$ must be one of $0$, $2u^2$
or $4u^2$.
The case $w=4u^2$ corresponds to $B_1+B_2+B_3+B_4=j$,
which clearly implies $\alpha=\beta=\gamma=\delta=y=u^2/2-u/2$.
In this case $(\diamondsuit)$ gives $q=0$.

The case $w=0$ corresponds to $B_1+B_2+B_3+B_4=0$,
which clearly implies $\alpha=\beta=\gamma=\delta=x=u^2/2-u$.
In this case $(\diamondsuit)$ gives $q=u^2/2-u$.

Finally, suppose $w=2u^2$.
Let $N_x$ be the number of $\alpha,\beta,\gamma,\delta$
that are equal to $x$, so $N_x\in \{0,1,2,3,4\}$.
Then $(\diamondsuit)$ gives $q=u^2/4-N_xu/4$, so each of the
five possibilities for $N_x$ gives the remaining
five possibilities for $q$.

\hfill $\Box$

\bigskip

\section{SDP Designs and Condition 2}

We now consider the question of which SDP designs 
$D(X,{\cal B})$ satisfy
condition 2 of Theorem \ref{bs1}, which states:
for any set ${\cal S}\subseteq {\cal B}$ 
of four blocks, an even number of the
four $3$-subsets of ${\cal S}$ have triple intersection size $x$.

Recall that any SDP design has parameters
$$(2^{2m},2^{2m-1}-2^{m-1},2^{2m-2}-2^{m-1})$$
and 2-rank $2m+2$,
and the derived design with respect to any block
has parameters
$$(2^{2m-1}-2^{m-1},2^{2m-2}-2^{m-1},2^{2m-2}-2^{m-1}-1) \quad  (\dag)$$
and 2-rank $2m+1$.
We will use the following result from McGuire and Ward \cite{MW}
(Corollary 4 and Theorem 8 there),
which characterises the derived designs of the
symplectic SDP designs by
their triple intersection sizes.

\begin{thm}\label{mw1} \cite{MW}
Let $D$ be a design with parameters $(\dag)$
and 2-rank $2m+1$.
Then $D$ is equivalent to
a derived design of the symplectic SDP design if and only
if all sizes of intersections of three blocks are
divisible by $2^{m-2}$.
\end{thm}

We now prove our main result.

\begin{thm}\label{mainthm}
Let $D(X,{\cal B})$  be a 
$(2^{2m},2^{2m-1}-2^{m-1},2^{2m-2}-2^{m-1})$ SDP design.
Then $D$ satisfies condition $2$ of Theorem $\ref{bs1}$
if and only if $D$ is equivalent to the symplectic SDP design.
\end{thm}

Proof:
We continue the notation of the proof of Theorem \ref{quadint}.
Condition 2 of Theorem \ref{bs1}
states that $N_x$ can only equal 0, 2, or 4, for all
choices of four blocks 
$B_1, B_2, B_3, B_4\in {\cal B}$.
By the proof of Theorem 3, 
this is equivalent to saying that only five
quadruple intersection sizes are allowed; the
two quadruple intersection sizes that are forbidden are (since
$u=2^{m-1}$)
$$u^2/4-u/4=2^{2m-4}-2^{m-3} \quad\textrm{and}\quad 
u^2/4-3u/4=2^{2m-4}-3\cdot 2^{m-3},$$ corresponding to
$N_x=1$ and $N_x=3$.

We further observe that the five allowable quadruple
intersection sizes are divisible by $2^{m-2}$, and the
two forbidden quadruple
intersection sizes are not divisible by $2^{m-2}$.
Therefore, $D$ will satisfy condition 2 of Theorem \ref{bs1}
if and only if all quadruple intersection sizes are
divisible by $2^{m-2}$.

Note that the  quadruple intersection sizes
of blocks in $D$ are the same as the 
triple intersection sizes of blocks in a derived
design of $D$.
Therefore, $D$ will satisfy condition 2 of Theorem \ref{bs1}
if and only if all triple intersection sizes of blocks in a derived
design of $D$ are
divisible by $2^{m-2}$.
Theorem \ref{mw1} implies that 
$D$ will satisfy condition 2 of Theorem \ref{bs1}
if and only if
any derived design of $D$
is equivalent to a derived design of the symplectic SDP design.
By a result of Jungnickel and Tonchev \cite{JT},
non-isomorphic SDP designs have non-isomorphic derived designs,
so a derived design of $D$ is equivalent to a derived design
of the symplectic design if and only if $D$ itself is equivalent 
to the symplectic design.

\hfill $\Box$

\bigskip

\section{Symplectic Designs and Condition 3}

We now show that the symplectic designs satisfy 
condition 3 of Theorem \ref{bs1}, which states:
there exists a 1-1 correspondence 
$\phi : X \longrightarrow {\cal B}$
with the property that for any three points
$a,b,c \in X$,
the number of blocks containing $\{a,b,c\}$ is
$|\phi(a)\cap \phi(b)\cap \phi(c)|$.

Recall that a \emph{polarity} of a square design $D(X,{\cal B})$ is
a bijection $\sigma : X \longrightarrow {\cal B}$
such that $\sigma \circ \sigma$ is the identity and
$p\in \sigma(q)$ if and only if $q\in \sigma(p)$,
for all $p, q \in X$.
A square design has a polarity if and only if
it has a symmetric incidence matrix, with respect
to some ordering of the points and blocks.

\begin{lemma}\label{polar}
Let $D(X,{\cal B})$  be a square design with a polarity.
Then there exists a 1-1 correspondence 
$\phi : X \longrightarrow {\cal B}$
with the property that for any three points
$a,b,c \in X$,
the number of blocks containing $\{a,b,c\}$ is
$|\phi(a)\cap \phi(b)\cap \phi(c)|$.
\end{lemma}

Proof: Let $\phi$ be the polarity of $D$.
Then, for $p, a, b, c, \in X$,
it follows from the definition of a polarity that
$$p\in (\phi(a) \cap \phi(b) \cap \phi(c)) \iff
\{a, b, c\} \subseteq \phi(p).$$

\hfill \hfill $\Box$

\begin{thm}\label{cond3}
The symplectic SDP designs satisfy condition 3 of
Theorem \ref{bs1}.
\end{thm}

The proof follows from Lemma \ref{polar} and the fact that
the symplectic SDP designs have a polarity
(see \cite{K1}, or \cite{CvL} page 78).

\bigskip

We now combine Theorems \ref{mainthm} and \ref{cond3}
to give our characterisation of the SDP designs
satisfying all of conditions 1, 2, and 3, of
Theorem \ref{bs1}.

\begin{thm}   Let $D$ be an SDP design.  Then $D$
satisfies  conditions $1, 2$ and $3$
of Theorem $1$ if and only if
$D$ is the symplectic SDP design.
\end{thm}

\bigskip

{\bf Acknowledgement}  We thank Wayne Broughton for
very helpful comments.

\newpage

Mailing address of contact author:

\bigskip

Gary McGuire

Department of Mathematics

NUI Maynooth

Co. Kildare

Ireland

\bigskip

Telephone: $+$ 353-1-708-3914

Fax: $+$ 353-1-708-3913

email:  gary.mcguire@may.ie

\end{document}